\documentclass{svmult}
\usepackage{amsmath,amssymb,bbm,graphicx}

\smartqed

\begin{document}

\title*{Long Memory in Nonlinear Processes}

  \author{Rohit Deo\inst{1} \and Mengchen Hsieh\inst{1} \and Clifford M.
    Hurvich\inst{1} \and Philippe Soulier\inst{2}}
  \institute{New York University, 44 W. 4'th Street, New York NY
    10012, USA \texttt{\{rdeo,mhsieh,churvich\}@stern.nyu.edu} 
        \and Universit\'e Paris X, 200 avenue de la R\'epublique, 92001
    Nanterre cedex, France \texttt{philippe.soulier@u-paris10.fr} 
}

\maketitle

\section{Introduction}

It is generally accepted that many time series of practical interest
exhibit strong dependence, i.e., long memory. For such series, the
sample autocorrelations decay slowly and log-log periodogram plots
indicate a straight-line relationship. This necessitates a class of
models for describing such behavior. A popular class of such models is
the autoregressive fractionally integrated moving average (ARFIMA)
(see \cite{adenstedt:1974}, \cite{granger:joyeux:1980}),
\cite{hosking:1981}, which is a linear process.  However, there is
also a need for nonlinear long memory models.  For example, series of
returns on financial assets typically tend to show zero correlation,
whereas their squares or absolute values exhibit long memory. See,
e.g., \cite{ding:granger:engle:1993}.  Furthermore, the search for a
realistic mechanism for generating long memory has led to the
development of other nonlinear long memory models. (Shot noise,
special cases of which are Parke, Taqqu-Levy, etc). In this chapter,
we will present several nonlinear long memory models, and discuss the
properties of the models, as well as associated parametric and
semiparametric estimators.

Long memory has no universally accepted definition; nevertheless, the most
commonly accepted definition of long memory for a weakly stationary
process $X=\{X_t,\ t\in\mathbb{Z}\}$ is the regular variation of the
autocovariance function: there exist $H\in(1/2,1)$ and a slowly
varying function $L$ such that
\begin{gather}  \label{eq:deflrd}
  \mathrm{cov}(X_0,X_t) = L(t) |t|^{2H-2} \; .
\end{gather}
Under this condition, it holds that:
\begin{gather}  \label{eq:varsum}
     \lim_{n\to\infty} n^{-2H} L(n)^{-1}   \mathrm{var} \left(
     \sum_{t=1}^n X_t \right) = 1/(2H(2H-1)).
\end{gather}
The condition~(\ref{eq:varsum}) does not imply~(\ref{eq:deflrd}).
Nevertheless, we will take~(\ref{eq:varsum}) as an alternate
definition of long memory.  In both cases, the index $H$ will be
referred to as the {\em Hurst index} of the process $X$.  This
definition can be expressed in terms of the parameter $d=H-1/2$, which
we will refer to as the {\em memory parameter}.  The most famous long
memory processes are fractional Gaussian noise and the $ARFIMA(p,d,q)$
process, whose memory parameter is $d$ and Hurst index is $H=1/2+d$. See
for instance \cite{taqqu:2003} for a definition of these processes.

The second-order properties of a stationary process are not sufficient
to characterize it, unless it is a Gaussian process. Processes which
are linear with respect to an i.i.d. sequence (strict sense linear
processes) are also relatively well characterized by their
second-order structure. In particular, weak convergence of the partial
sum process of a Gaussian or strict sense linear long memory processes
$\{X_t\}$ with Hurst index $H$ can be easily derived. Define $S_n(t) =
\sum_{k=1}^{[nt]} (X_k-\mathbb{E}[X_k])$ in discrete time or $S_n(t) =
\int_0^{nt} (X_s-\mathbb{E}[X_s]) \D s$ in continuous time. Then
$\mathrm{var}(S_n(1))^{-1/2} S_n(t)$ converges in distribution to a
constant times the fractional Brownian motion with Hurst index $H$,
that is the Gaussian process $B_H$ with covariance function
\begin{gather*}
  \mathrm{cov}(B_H(s),B_H(t)) = \frac12 \{|s|^{2H} - |t-s|^{2H} +
  t^{2H}\} \; .
\end{gather*}

In this paper, we will introduce nonlinear long memory processes,
whose second order structure is similar to that of Gaussian or linear
processes, but which may differ greatly from these processes in many
other aspects. In Section~\ref{sec:models}, we will present these
models and their second-order properties, and the weak convergence of
their partial sum process. These models include conditionally
heteroscedastic processes (Section~\ref{sec:bilinear}) and models
related to point processes (Section~\ref{sec:renewal}). In
Section~\ref{sec:estimation}, we will consider the problem of
estimating the Hurst index or memory parameter of these processes.

\section{Models}
\label{sec:models}

\subsection{Conditionally heteroscedastic  models}
\label{sec:bilinear}
These models are defined by 
\begin{gather}
  X_t=\sigma _t v_t \; , \label{eq:chm}
\end{gather} 
where $\{v_t\}$ is an independent identically distributed series with
finite variance and $\sigma_t^2$ is the so-called volatility.  
%% If $\sigma_t^2$ is measurable with respect to the past of the process,
%% then $\sigma_t^2 = \mathbb{E}[X_t^2 \mid X_{t-1}, X_{t-2},\dots]$ is
%% the conditional variance. 
We now give examples.
\paragraph{LMSV and LMSD}
The Long Memory Stochastic Volatility (LMSV) and Long Memory
Stochastic Duration (LMSD) models are defined by
Equation~(\ref{eq:chm}), where $\sigma _t^2=\exp({h_t})$ and $\{h_t\}$
is an unobservable Gaussian long memory process with memory parameter
$d\in(0,1/2)$, independent of $\{v_t\}$.  The multiplicative
innovation series $\{v_t\}$ is assumed to have zero mean in the LMSV
model, and positive support with unit mean in the LMSD model. The LMSV
model was first introduced by \cite{breidt:crato:delima:1998} and
\cite{harvey:1998} to describe returns on financial assets, while the LMSD
model was proposed by \cite{deo:hsieh:hurvich:2005} to describe
durations between transactions on stocks.

Using the moment generating function of a Gaussian distribution, it
can be shown (see \cite{harvey:1998}) for the LMSV/LMSD model that for
any real $s$ such that $\mathbb{E}[|v_t|^s] < \infty$,
\[
\rho_s(j) \sim C_s j^{2d-1} \qquad j \rightarrow \infty ,
\]
where $\rho_s(j)$ denotes the autocorrelation of $\left\{ \left|
    x_t\right| ^s\right\} $ at lag $j$, with the convention that $s=0$
corresponds to the logarithmic transformation.  As shown in
\cite{surgailis:viano:2002}, the same result holds under more general
conditions without the requirement that $\{h_t\}$ be Gaussian.

In the LMSV model, assuming that $\{h_t\}$ and $\{v_t\}$ are functions
of a multivariate Gaussian process, \cite{robinson:2001} obtained
similar results on the autocorrelations of $\{|X_t|^s\}$ with $s>0$
even if $\{h_t\}$ is not independent of $\{v_t\}$. Similar results
were obtained in \cite{surgailis:viano:2002}, allowing for dependence
between $\{h_t\}$ and~$\{v_t\}$.

The LMSV process is an uncorrelated sequence, but powers of LMSV or
LMSD may exhibit long memory.  \cite{surgailis:viano:2002} proved the
convergence of the centered and renormalized partial sums of any
absolute power of these processes to fractional Brownian motion with
Hurst index 1/2 in the case where they have short memory.

\paragraph{FIEGARCH}
The weakly stationary FIEGARCH model was proposed by
\cite{bollerslev:mikkelsen:1996}. The FIEGARCH model, which is
observation-driven, is a long-memory extension of the EGARCH
(exponential GARCH) model of \cite{nelson:1991}. The FIEGARCH model
for returns $\{X_t\}$ takes the form~\ref{sec:bilinear} innovation
series $\{v_t\}$ are i.i.d. with zero mean and a symmetric
distribution, and
\begin{equation}\label{eq.fiegarch1}
\log \sigma_t^2 = \omega + \sum_{j=1}^\infty a_j g(v_{t-j})
\end{equation}
with $g(x)=\theta x + \gamma (|x|-\mathbb{E} |v_t|)$, $\omega >0$,
$\theta\in\mathbb{R}$, $\gamma\in\mathbb{R}$, and real constants $a_j$
such that the process $\{\log \sigma_t^2\}$ has long memory with
memory parameter $d\in(0,1/2)$.  If $\theta$ is nonzero, the model
allows for a so-called leverage effect, whereby the sign of the
current return may have some bearing on the future volatility. In the
original formulation of \cite{bollerslev:mikkelsen:1996}, the
$\{a_j\}$ are the $AR(\infty)$ coefficients of an $ARFIMA(p,d,q)$
process.

As was the case for the LMSV model, here we can once again express the
log squared returns as in (\ref{eq:SigNoise}) with $\mu = \mathbb{E}
[\log v_t^2]+\omega$, $u_t = \log v_t^2 - \mathbb{E} [ \log v_t^2 ]$,
and $h_t=\log \sigma_t^2 - \omega$. Here, however, the processes
$\{h_t\}$ and $\{u_t\}$ are not mutually independent. The results of
\cite{surgailis:viano:2002} also apply here, and in particular, the
processes $\{|X_t|^u\}$, $\{\log(X_t^2)\}$ and $\{\sigma_t\}$
have the same memory parameter $d$.

\paragraph{ARCH($\infty$) and FIGARCH}

In ARCH($\infty$) models, the innovation series $\{v_t\}$ is assumed
to have zero mean and unit variance, and the conditional variance is
taken to be a weighted sum of present and past squared returns:
\begin{gather} \label{eq:volatilityarchinfty}
  \sigma_t^2 = \omega  + \sum_{k=1}^\infty a_j X_{t-j}^2 \; ,
\end{gather}
where $\omega, a_j, j=1,2,\dots$ are nonnegative constants. The
general framework leading to~(\ref{eq:chm})
and~(\ref{eq:volatilityarchinfty}) was introduced by
\cite{robinson:1991}.  \cite{kazakevicius:leipus:2003} have shown that
$\sum_{j=1}^\infty a_j \leq 1$ is a necessary condition for existence
of a strictly stationary solution to
equations~(\ref{eq:chm}),~(\ref{eq:volatilityarchinfty}), while
\cite{giraitis:kokoszka:leipus:2000} showed that $\sum_{j=1}^\infty
a_j < 1$ is a sufficient condition for the existence of a strictly
stationary solution. If $\sum_{j=1}^\infty a_j = 1$, the existence of
a strictly stationary solution has ben proved by
\cite{kazakevicius:leipus:2003} only in the case where the
coefficients $a_j$ decay exponentially fast. In any case, if a
stationary solution exists, its variance, if finite, must be equal to
$\omega(1-\sum_{k=1}^\infty a_k)^{-1}$, so that it cannot be finite if
$\sum_{k=1}^\infty a_k=1$ and $\omega>0$. If $\omega=0$, then the
process which is identically equal to zero is a solution, but it is
not known whether a nontrivial solution exists.

In spite of a huge literature on the subject, the existence of a
strictly or weakly stationary solution
to~(\ref{eq:chm}),~(\ref{eq:volatilityarchinfty}) such that
$\{\sigma_t^2\}$, $\{|X_t|^u\}$ or $\{\log(X_t^2)\}$ has long memory
is still an open question. If $\sum_{j=1}^\infty a_j<1$, and the
coefficients $a_j$ decay sufficiently slowly,
\cite{giraitis:kokoszka:leipus:2000} found that it is possible in such
a model to get hyperbolic decay in the autocorrelations $\{\rho_r\}$
of the squares, though the rates of decay they were able to obtain
were proportional to $r^{-\theta}$ with $\theta > 1$. Such
autocorrelations are summable, unlike the autocorrelations of a
long-memory process with positive memory parameter. For instance, if
the weights $\{a_j\}$ are proportional to those given by the
$AR(\infty)$ representation of an ARFIMA($p,d,q$) model, then $\theta
= -1-d$.  If $\sum_{j=1}^\infty a_j = 1$, then the process has
infinite variance so long memory as defined here is irrelevant.

Let us mention for historical interest the FIGARCH (fractionally
integrated GARCH) model which appeared first in
\cite{baillie:bollerslev:mikkelsen:1996}. In the FIGARCH model, the
weights $\{a_j\}$ are given by the $AR(\infty)$ representation of an
ARFIMA($p,d,q$) model, with $d \in (0,1/2)$, which implies that
$\sum_{j=1}^\infty a_j = 1$, hence the very existence of FIGARCH
series is an open question, and in any case, if it exists, it cannot
be weakly stationary. The lack of weak stationarity of the FIGARCH
model was pointed out by \cite{baillie:bollerslev:mikkelsen:1996}.
Once again, at the time of writing this paper, we are not aware of any
rigorous result on this process or on any ARCH($\infty)$ process with
long memory.

%% Since long memory is typically defined in terms of the
%% autocorrelations of the series, one cannot properly say that the
%% $FIGARCH$ model has long memory (either in levels or in squares), as
%% its autocorrelations are not well defined.

%% Zaffaroni (2000) considered the possibility of taking the constant
%% $\omega$ to be zero, but still was unable to establish whether it is
%% possible to get long memory in the squares of such a process. Davidson
%% (2004) has summarized much of the work on this question, and states
%% that one should look beyond the FIGARCH class if one wants to have
%% long memory in volatility.

\paragraph{LARCH}
Since the ARCH structure (appearently) fails to produce long memory,
an alternative definition of heteroskedasticity has been considered in
which long memory can be proved rigorously.
\cite{giraitis:surgailis:2002} considered models which satisfy the
equation $X_t = \zeta_t A_t + B_t$, where $\{\zeta_t\}$ is a sequence
of i.i.d. centered random variables with unit variance and $A_t$ and
$B_t$ are linear in $\{X_t\}$ instead of quadratic as in the ARCH
specification. This model nests the LARCH model introduced by
\cite{robinson:1991}, obtained for $B_t \equiv 0$. The advantage of
this model is that it can exhibit long memory in the conditional mean
$B_t$ and/or in the conditional variance $A_t$, possibly with
different memory parameters.  See \cite[Corollary
4.4]{giraitis:surgailis:2002}.  The process $\{X_t\}$ also exhibits
long memory with a memory parameter depending on the memory parameters
of the mean and the conditional variance \cite[Theorem
5.4]{giraitis:surgailis:2002}. If the conditional mean exhibits long
memory, then the partial sum process converges to the fractional
Brownian motion, and it converges to the standard Brownian motion
otherwise. See \cite[Theorem 6.2]{giraitis:surgailis:2002}.  The
squares $\{X_t^2\}$ may also exhibit long memory, and their partial
sum process converge either to the fractional Brownian motion or to a
non Gaussian self-similar process.  This family of processes is thus
very flexible. An extension to the multivariate case is given in
\cite{doukhan:teyssiere:winant:2005}.

\vspace{2em}

We conclude this section by the following remark. Even though these
processes are very different from Gaussian or linear processes, they
share with weakly dependent processes the Gaussian limit and the fact
that weak limits and $L^2$ limits have consistent normalisations, in
the sense that, if $\xi_n$ denotes one of the usual statistics
computed on a time series, there exists a sequence $v_n$ such that
$v_n \xi_n$ converges weakly to a non degenerate distribution and
$v_n^2 \mathbb{E}[\xi_n^2]$ converges to a positive limit (which is the
variance of the asymptotic distribution). In the next subsection, we
introduce models for which this is no longer true.

\subsection{Shot noise processes}
\label{sec:renewal}
 
General forms of the shot-noise process have been considered for a long
time; see for instance \cite{takacs:1954}, \cite{daley:1971}. Long
memory shot noise processes have been introduced more recently; an
early reference seems to be \cite{giraitis:molchanov:surgailis:1993}.
We present some examples of processes related to shot noise which may
exhibit long memory. For simplicity and brevity, we consider only
stationary processes.

Let $\{t_j, \; j \in \mathbb{Z}\}$ be the points of a stationary point
process on the line, numbered for instance in such a way that $t_{-1}
< 0 \leq t_0$, and for $t\geq0$, let
$N(t)=\sum_{j\geq0}\mathbbm{1}_{\{t_j\leq t\}}$ be the number of points
between time zero and $t$. Define then
\begin{gather} \label{eq:defshotnoise}
  X_t = \sum_{j \in \mathbb{Z}} \epsilon_j \mathbbm{1}_{\{ t_j \leq t < t_j +
    \eta_j\}}, \ \ \ t \geq 0.
\end{gather} 
In this model, the shocks $\{\epsilon_j\}$ are an i.i.d. sequence;
they are generated at birth times $\{t_j\}$ and have durations
$\{\eta_j\}$.
The observation at time $t$ is the sum of all surviving
present and past shocks. In model~(\ref{eq:defshotnoise}), we can take
time to be continuous, $t\in\mathbb{R}$ or discrete, $t\in\mathbb{Z}$. This will
be made precise later for each model considered. We now describe
several well known special cases of model~(\ref{eq:defshotnoise}).

\begin{enumerate}
\item Renewal-reward process; \cite{taqqu:levy:1986},  \cite{liu:2000}. \\
  The durations are exactly the interarrival times of the renewal
  process: $\eta_0 = t_0$, $\eta_j = t_{j+1} - t_j$, and the shocks
  are independent of their birth times.  Then there is exactly one
  surviving shock at time~$t$:
  \begin{gather} \label{eq:defrenewalreward}
    X_t = \epsilon_{N(t)}.
  \end{gather}
\item ON-OFF model; \cite{taqqu:willinger:sherman:1997}. \\
  This process consists of alternating ON and OFF periods with
  independent durations. Let $\{\eta_k\}_{\geq1}$ and
  $\{\zeta_k\}_{k\geq1}$ be two independent i.i.d. sequences of
  positive random variables with finite mean. Let $t_0$ be independent
  of these sequences and define $t_j = t_0 + \sum_{k=1}^j (\eta_k +
  \zeta_k)$. The shocks $\epsilon_j$ are deterministic and equal to 1.
  Their duration is $\eta_j$. The $\eta_j$s are the ON periods and the
  $\zeta_j$s are the OFF periods.  The first interval $t_0$ can also be
  split into two successive ON and OFF periods $\eta_0$ and $\zeta_0$.
  The process $X$ can be expressed as
  \begin{gather} \label{eq:defonoff}
    X_t = \mathbbm{1}_{\{t_{N(t)} \leq t < t_{N(t)}+\eta_{N(t)}\}}.
  \end{gather}

\item Error duration process; \cite{parke:1999}. \\
  This process was introduced to model some macroeconomic data. The
  birth times are deterministic, namely $t_j = j$, the durations
  $\{\eta_j\}$ are i.i.d. with finite mean and
  \begin{gather} \label{eq:defparkec}
    X_t = \sum_{j \leq t} \epsilon_j \mathbbm{1}_{\{t < j + \eta_j\}}. 
  \end{gather}
%  (***You should state here your assumptions on the $\eta_j$***)
%  
% no, I do not want to say more than they are i.i.d. here. 
% here, the model is just defined, the precise assumptions 
% come in the next paragraph; 
%
\item Infinite Source Poisson model.\\
  If the $t_j$ are the points of a homogeneous Poisson process, the
  durations $\{\eta_j\}$ are i.i.d. with finite mean and
  $\epsilon_j\equiv 1$, we obtain the infinite source Poisson model or
  M/G/$\infty$ input model considered among others in
  \cite{mikosch:resnick:rootzen:stegeman:2002}.
%%(***You should state   here your assumptions on the $\eta_j$***)
%%
%% same answer
%%

  \cite{maulik:resnick:rootzen:2002} have considered a variant of this
  process where the shocks (referred to as transmission rates in this
  context) are random, and possibly contemporaneously dependent with
  durations.

  \end{enumerate}  
  In the first two models, the durations satisfy $\eta_j \leq
  t_{j+1}-t_j$, hence are not independent of the point process of
  arrivals (which is here a renewal process). Nevertheless $\eta_j$ is
  independent of the past points $\{t_k, \; k \leq j \}$. The process
  can be defined for all $t\geq0$ without considering negative birth
  times and shocks. In the last two models, the shocks and durations
  are independent of the renewal process, and any past shock may
  contribute to the value of the process at time $t$.

\paragraph{Stationarity and second order properties} 
$\bullet$ The renewal-reward process~(\ref{eq:defrenewalreward}) is strictly
stationary since the renewal process is stationary and the shocks are
i.i.d. It is moroever weakly stationary if the shocks have finite
variance. Then  $\mathbb{E}[X_t] = \mathbb{E}[\epsilon_1]$ and 
\begin{gather} \label{eq:covrr}
  \mathrm{cov}(X_0,X_t) = \mathbb{E}[\epsilon^2] \; \mathbb{P}(\eta_0>t) = \lambda
  \mathbb{E}[\epsilon_1^2] \; \mathbb{E}[(\eta_1-t)_+]\; ,
\end{gather}
where $\eta_0$ is the delay distribution and $\lambda =
\mathbb{E}[(t_1-t_0)]^{-1}$ is intensity of the stationary renewal process.
Note that this relation would be true for a general stationary point
process.  Cf. for instance \cite{taqqu:levy:1986} or
\cite{hsieh:hurvich:soulier:2003}.\\

\noindent $\bullet$ The stationary version of the ON-OFF was studied in
\cite{heath:resnick:samorodnitsky:1998}. The first On and OFF period
$\eta_0$ and $\zeta_0$ can be defined in such a way that the process
$X$ is stationary.  Let $F_{\mathrm{on}}$ and $F_{\mathrm{off}}$ be
the distribution functions of the ON and OFF periods $\eta_1$ and
$\zeta_1$. \cite[Theorem 4.3]{heath:resnick:samorodnitsky:1998} show that if
$1-F_{\mathrm{on}}$ is regularly varying with index $\alpha\in(1,2)$
and $1-F_{\mathrm{off}}(t) = o(F_{\mathrm{on}}(t))$ as $t\to\infty$, then
\begin{gather} \label{eq:covonoff}
 \mathrm{cov}(X_0,X_t) \sim c  \mathbb{P}(\eta_0>t) = \; c \lambda \mathbb{E}[(\eta_1-t)_+]\; ,
\end{gather}

\noindent $\bullet$ Consider now the case when the durations are independent of the birth
times. To be precise, assume that $\{(\eta_j,\epsilon_j)\}$ is an
i.i.d. sequence of random vectors, independent of the stationary point
process of points $\{t_j\}$. Then the process $\{X_t\}$ is strictly
stationary as long as $\mathbb{E}[\eta_1] < \infty$, and has finite variance
if $\mathbb{E}[\epsilon_1^2 \eta_1] < \infty$. Then $\mathbb{E}[X_t] = \lambda
\mathbb{E}[\epsilon_1\eta_1]$ and 
\begin{align*}
  \mathrm{cov}(X_0,X_t) & = \lambda \,  \mathbb{E}[\epsilon_1^2 \, (\eta_1-t)_+] \\
  + & \; \{ \mathrm{cov}(\epsilon_1 \, N(-\eta_1,0], \epsilon_2 \,
  N(t-\eta_2,t]) - \lambda \mathbb{E}[\epsilon_1 \epsilon_2 \, (\eta_1
  \wedge (\eta_2-t)_+]\} \; ,
\end{align*}
where $\lambda$ is the intensity of the stationary point process, i.e.
$\lambda^{-1} = \mathbb{E}[t_0]$.  The last term has no known general
expression for a general point process, but it vanishes in two
particular cases:
\begin{itemize} 
\item[-] if $N$ is a homogeneous Poisson point process;
\item[-] if $\epsilon_1$ is centered and independent of $\eta_1$.
\end{itemize}
In the latter case~(\ref{eq:covrr}) holds, and in the former case, we
obtain a formula which generalizes~(\ref{eq:covrr}):
\begin{gather} \label{eq:covgeneral}
  \mathrm{cov}(X_0,X_t)  = \lambda \,  \mathbb{E}[\epsilon_1^2 \, (\eta_1-t)_+] \; .
\end{gather}
We now see that second order long memory can be obtained
if~(\ref{eq:covrr}) holds and the durations have regularly varying
tails with index $\alpha \in (1,2)$ or, 
\begin{gather} \label{eq:nonindep}
  \mathbb{E}[\epsilon_1^2 \mathbbm{1}_{\{\eta_1 > t \}} ] = \ell(t) t^{-\alpha} \; .
\end{gather}
Thus, if~(\ref{eq:nonindep}) and either~(\ref{eq:covonoff})
or~(\ref{eq:covgeneral}) hold, then $X$ has long memory with Hurst
index $H = (3-\alpha)/2$ since
\begin{gather} \label{eq:cov}
  \mathrm{cov}(X_0,X_t) \sim  \, \frac{\lambda}{\alpha - 1} \; \ell(t)
  t^{1-\alpha} \; .
\end{gather}
Examples of interest in teletraffic modeling where $\epsilon_1$ and
$\eta_1$ are not independent but~(\ref{eq:nonindep}) holds are provided
in \cite{maulik:resnick:rootzen:2002} and
\cite{fay:roueff:soulier:2005}.

We conjecture that~(\ref{eq:cov}) holds in a more general framework,
at least if the interarrival times of the point process have finite
variance.

%%    \begin{itemize}
%%    \item $X$ is the renewal-reward process with regularly varying
%%      interarrival times with index $\alpha\in(1,2)$;
%%    \item $X$ is the ON-OFF process and the ON duration is regularly
%%      varying with index $\alpha\in(1,2)$;
%%    \item The shocks and durations are i.i.d. and indepedent of the
%%      stationary point process $\{t_j\}$; the shocks are independent of
%%      the durations which are regularly varying with index
%%      $\alpha\in(1,2)$; 
%%    \item The shocks and durations are i.i.d. and indepedent of the
%%      homogeneous Poisson point process $\{t_j\}$; the joint
%%      distribution of contemporaneous shocks and durations
%%      satisfy~(\ref{eq:nonindep}).
%%    \end{itemize}

\paragraph{Weak convergence of partial sums}
This class of long memory process exhibits a very distinguishing
feature. Instead of converging weakly to a process with finite
variance, dependent stationary increments such as the fractional
Brownian motion, the partial sums of some of these processes have been
shown to converge to an $\alpha$-stable Levy process, that is, an
$\alpha$-stable process with independent and stationary increment.
Here again there is no general result, but such a convergence is easy
to prove under restrictive assumptions. Define 
\begin{gather*}
  S_T(t) =  \int_0^{Tt} \{ X_s - \mathbb{E}[X_s] \} \, \D s \; .
\end{gather*}
Then it is known in the particular cases described above that the
finite dimensional distributions of the process $\ell(T)
T^{-1/\alpha}S_T$ (for some slowly varying function $\ell$) converge
weakly to those of an $\alpha$-stable process. This was proved in
\cite{taqqu:levy:1986} for the renewal reward process, in
\cite{mikosch:resnick:rootzen:stegeman:2002} for the ON-OFF and
infinite source Poisson processes when the shocks are constant.  A
particular case of dependent shocks and durations is considered in
\cite{maulik:resnick:rootzen:2002}.  \cite{hsieh:hurvich:soulier:2003}
proved the result in discrete time for the error duration process; the
adaptation to the continuous time framework is straightforward.  It is
also probable that such a convergence holds when the underlying point
process is more general.

Thus, these processes are examples of second order long memory process
with Hurst index $H\in(1/2,1)$ such that $T^{-H}S_T(t)$ converges in
probability to zero. This behaviour is very surprising and might be
problematic in statistical applications, as illustrated in
Section~\ref{sec:estimation}.

It must also be noted that convergence does not hold in the space
$\mathcal{D}$ of right-continuous, left-limited functions endowed with
the $J_1$ topology, since a sequence of processes with continuous path
which converge in distribution in this sense must converge to a
process with continuous paths.  It was proved in \cite[Theorem
4.1]{resnick:vandenberg:2000} that this convergence holds in the $M_1$
topology for the infinite source Poisson process. For a definition and
application of the $M_1$ topology in queuing theory, see
\cite{whitt:2002}.

\paragraph{Slow growth and fast growth}
Another striking feature of these processes is the slow growth versus
fast growth phenomenon, first noticed by \cite{taqqu:levy:1986} for
the renewal-rewrd process and more rigorously investigated by
\cite{mikosch:resnick:rootzen:stegeman:2002} for the ON-OFF and
infinite source Poisson process\footnote{Actually, in the case of the
  Infinite Source Poisson process,
  \cite{mikosch:resnick:rootzen:stegeman:2002} consider a single
  process but with an increasing rate $\lambda$ depending on $T$,
  rather than superposition of independent copies. The results
  obtained are nevertheless of the same nature.}. Consider $M$
independent copies $X^{(i)}$, $1, \leq i \leq M$ of these processes
and denote
\begin{gather*}
  A_{M,T}(t) = \sum_{i=1}^M \int_0^{Tt} \{ X_s^{(i)} - \mathbb{E}[X_s] \} \,
  \D s \; .
\end{gather*}
If $M$ depends on $T$, then, according to the rate growth of $M$ with
respect to $T$, a stable or Gaussian limit can be obtained. More
precisely, the slow growth and fast growth conditions are, up to
slowly varying functions $M T^{1-\alpha} \to 0$ and $MT^{1-\alpha} \to
\infty$, respectively. In other terms, the slow and fast growth
conditions are characterized by $\mathrm{var}(A_{M,T}(1)) \ll b(MT)$ and
$\mathrm{var}(A_{M,T}(1)) \gg b(MT)$, respectively, where $b$ is the inverse
of the quantile function of the durations.

Under the slow growth condition, the
finite dimensional distributions of $L(MT) (MT)^{-1/\alpha} A_{M,T}$
converge to those of a Levy $\alpha$-stable process, where $L$ is a
slowly varying function. Under the fast growth condition, the sequence
of processes $T^{-H} \ell^{-1/2}(T) M^{-1/2} A_{M,T}$ converges, in
the space $\mathcal{D}(\mathbb{R}_+)$ endowed with the $J_1$ topology, to
the fractional Brownian motion with Hurst index $H = (3-\alpha)/2$.
It is thus seen that under the fast growth condition, the behaviour of
a Gaussian long memory process with Hurst index $H$ is recovered.

\paragraph{Non stationary versions}

If the sum defining the process $X$ in~(\ref{eq:defshotnoise}) is
limited to non negative indices $j$, then the sum has always a finite
number of terms and there is no restriction on the distribution of the
interarrival times $t_{j+1}-t_j$ and the durations $\eta_j$. These
models can then be nonstationary in two ways: either because of
initialisation, in which case a suitable choice of the initial
distribution can make the process stationary; or because these
processes are non stable and have no stationary distribution. The
latter case arises when the interarrival times and/or the durations
have infinite mean. These models were studied by
\cite{resnick:rootzen:2000} and \cite{mikosch:resnick:2004} in the
case where the point process of arrivals is a renewal process.
contrary to the stationry case, where heavy tailed durations imply non
Gaussian limits, the limiting process of the partial sums has non
stationary increments and can be Gaussian in some cases.

\subsection{Long Memory in Counts}

The time series of counts of the number of transactions in a given
fixed interval of time is of interest in financial econometrics.
Empirical work suggests that such series may possess long memory.
See \cite{deo:hsieh:hurvich:2005}. Since the counts are induced by
the durations between transactions, it is of interest to study the
properties of durations, how these properties generate long memory
in counts, and whether there is a connection between potential
long memory in durations and long memory in counts.

The event times determine a counting process $N(t)=$ Number of
events in $(0,t]$. Given any fixed clock-time spacing $\Delta t
>0$, we can form the time series $\{\Delta
N_{t^\prime}\}=\{N(t^\prime \Delta t)-N[(t^\prime-1)\Delta t]\}$
for $t^\prime = 1,2,\ldots$, which counts the number of events in
the corresponding clock-time intervals of width $\Delta t$. We
will refer to the $\{\Delta N_{t^\prime}\}$ as the $counts$. Let
$\tau_k >0$ denote the waiting time (duration) between the
$k-1$'st and the $k$'th transaction.

We give some preliminary definitions taken from
\cite{daley:vere-jones:2003}.
\begin{definition}

A point process $N(t)=N(0,t]$ is \textit{stationary} if for every
$r=1,2,\ldots$ and all bounded Borel sets $A_1,\ldots ,A_r$, the
joint distribution of $\{N(A_1+t) ,\ldots , N(A_r+t)\}$ does not
depend on $t\in [0,\infty)$.

A second order stationary point process is \textit{long-range
count dependent} ($LRcD$) if 
$$
\lim_{t\rightarrow \infty} \frac{\mathrm{var}(N(t))}{t} = \infty \;
.
$$

A second order stationary point process $N(t)$ which is $LRcD$ has
\textit{Hurst index} $H \in (1/2,1)$ given by
\[
H = \sup \{ h: \limsup _{t \rightarrow \infty} \frac {\mathrm{var} (N(t))}
{t^{2h}} = \infty \} \; .
\]
\end{definition}

Thus if the counts $\{\Delta
N_{t^\prime}\}_{t^\prime=-\infty}^\infty$ on intervals of any
fixed width $\Delta t>0$ are LRD with memory parameter $d$ then
the counting process $N(t)$ must be LRcD with Hurst index
$H=d+1/2$. Conversely, if $N(t)$ is an LRcD process with Hurst
index $H$, then $\{\Delta N_{t^\prime}\}$ cannot have
exponentially decaying autocorrelations, and under the additional
assumption of a power law decay of these autocorrelations,
$\{\Delta N_{t^\prime}\}$ is LRD with memory parameter $d=H-1/2$.

There exists a probability measure $P^0$ under which the doubly
infinite sequence of durations $\{\tau_k\}_{k=-\infty}^{\infty}$ are a
stationary time series, i.e., the joint distribution of any
subcollection of the $\{\tau_k\}$ depends only on the lags between the
entries. On the other hand, the point process $N$ on the real line is
stationary under the measure $P$. A fundamental fact about point
processes is that in general (a notable exception is the Poisson
process) there is no single measure under which both the point process
$N$ and the durations $\{\tau_k\}$ are stationary, i.e., in general
$P$ and $P^0$ are not the same. Nevertheless, there is a one-to-one
correspondence between the class of measures $P^0$ that determine a
stationary duration sequence and the class of measures $P$ that
determine a stationary point process. The measure $P^0$ corresponding
to $P$ is called the \textit{Palm distribution}. The counts are
stationary under $P$, while the durations are stationary under $P^0$.

We now present an important theoretical result obtained by
\cite{daley:1999}.

\begin{theorem} \label{theo:first}
  A stationary \textbf{renewal} point process is LRcD and has Hurst
  index $H=(1/2)(3-\alpha)$ under $P$ if the interarrival time has
  tail index $ 1< \alpha < 2$ under $P^0$.
%% $H = \frac{1}{2}(3-\alpha)$ if the interarrival time has
%% tail index $1<\alpha<2$ .
\end{theorem}

Theorem \ref{theo:first} establishes a connection between the tail
index of a duration process and the persistence of the counting
process. According to the theorem, the counting process will be LRcD
if the duration process is $iid$ with infinite variance.  Here, the
memory parameter of the counts is completely determined by the tail
index of the durations.

This prompts the question as to whether long memory in the counts can
be generated solely by dependence in finite-variance durations. An
answer in the affirmative was given by
\cite{daley:rolski:vesilo:2000}, who provide an example outside of the
framework of the popular econometric models. We now present a theorem
on the long-memory properties of counts generated by durations
following the LMSD model. The theorem is a special case of a result
proved in \cite{deo:hurvich:soulier:wang:2005}, who give sufficient
conditions on durations to imply long memory in counts.

\begin{theorem} \label{theo:LMSD}
  If the durations $\{ \tau_k \}$ are generated by the LMSD process
  with memory parameter $d$, then the induced counting process $N(t)$
  has Hurst index $H=1/2+d$, i.e.  satisfies $\textrm{var}(N(t)) \sim
  Ct^{2d+1}$ under $P$ as $t \rightarrow \infty$ where $C>0$.
\end{theorem}

%% \subsection{Dynamical systems}
%% \label{sec:dynamic}

%% !!!!!!!!!!!!!!!!!!!!!!!   TO BE DONE? MAYBE

\section{Estimation of the Hurst index or memory parameter}
\label{sec:estimation}

A weakly stationary process with autocovariance function satisfying
(\ref{eq:deflrd}) has a spectral density $f$ defined by
\begin{gather} \label{eq:defspecdens}
f(x) = \frac1{2\pi} \sum_{t\in \mathbb{Z}} \gamma(t) \E^{\I tx} \; .
\end{gather}
This series converges uniformly on the compact subsets of
$[-\pi,\pi]\setminus\{0\}$ and in $L^1([-\pi,\pi],dx)$. Under some
strengthening of condition~(\ref{eq:deflrd}), the behaviour of the
function $f$ at zero is related to the rate of decay of~$\gamma$. For
instance, if we assume in addition that $L$ is ultimately monotone, we
obtain the following Tauberian result \cite[Proposition
4.1]{taqqu:2003}, with $d=H-1/2$.
\begin{gather} \label{eq:spectraldensity}
  \lim_{x \to 0} L(x)^{-1} x^{2d} f(x) = \pi^{-1} \Gamma(2d)
  \cos(\pi d).
\end{gather}
Thus, a natural idea is to estimate the spectral density in order to
estimate the memory paramter $d$. The statistical tools are the
discrete Fourier transform (DFT) and the periodogram, defined for a
sample $U_1,\dots,U_{n}$, as
\[
J_{n,j}^U = (2\pi n)^{-1/2} \sum_{t=1}^n U_t \E^{\I t w_j}, \ \ 
I_U(\omega_j) = |J_{n,j}^U|^2,
\]
where $\omega_j = 2j\pi/n$, $1 \leq j < n/2$ are the so-called Fourier
frequencies. (Note that for clarity the index $n$ is omitted from the
notation).  In the classical weakly stationary short memory case (when
the autocovariance function is absolutely summable), it is well known
that the periodogram is an asymptotically unbiased estimator of the
spectral density $f_U$ defined in~(\ref{eq:defspecdens}).  This is no
longer true for second order long memory processes.
\cite{hurvich:beltrao:1993} showed (in the case where the function $L$
is continuous at zero but the extension is straightforward) that for
any fixed positive integer $j$, there exists a positive constant
$c(k,H)$ such that
\[
\lim_{n\to\infty} \mathbb{E}[I_U(\omega_j)/f_U(\omega_j)] = c(j,H).
\]
The previous results are true for any second order long memory
process. Nevertheless, spectral method of estimation of the Hurst
parameter, based on the heuristic (but incorrect) assumption that
the renormalised DFTs $f_U^{-1/2}(\omega_j) J_{n,j}^U$ are i.i.d.
standard complex Gaussian have been proposed and theoretically
justifed in some cases. The most well known is the GPH estimator of
the Hurst index, introduced by \cite{geweke:porter-hudak:1983} and
proved consistent and asymptotically Gaussian for Gaussian long memory
processes by \cite{robinson:1995l} and for a restricted class of
linear processes by \cite{velasco:2000}. Another estimator, often
referred to as the local Whittle or GSE estimator was introduced by
\cite{kunsch:1987} and again proved consistent asymptotically Gaussian
by \cite{robinson:1995g} for linear long memory processes.

These estimators are built on the $m$ first log-periodogram ordinates,
%% and the GPH estimator boils down to a regression on $\log(j)$ or
%% $\log(|\sin(\omega_j/2)|)$, $1 \leq j \leq m$, 
where $m$ is an intermediate sequence, i.e.  $1/m+m/n \to 0$ as
$n\to\infty$. The choice of $m$ is irrelevant to consistency of the
estimator but has an influence on the bias.  The rate of
convergence of these estimators, when known, is typically slower than
$\sqrt n$. Trimming of the lowest frequencies, which means taking the
$\ell$ first frequencies out is sometimes used, but there is no
theoretical need for this practice, at least in the Gaussian case.
See \cite{hurvich:deo:brodsky:1998}. For nonlinear series, we are not
sure yet if trimming may be needed in general.

In the following subsections, we review what is known, both
theoretically and empirically, about these and related methods for the
different types of nonlinear processes described previsoulsy.

We start by describing the behaviour of the renormalized DFTs at low
frequencies, that is, when the index $j$ of the frequency $\omega_j$
remains fixed as $n\to\infty$.

\subsection{Low-Frequency DFTs of Counts from Infinite-Variance Durations}

To the best of our knowledge there is no model in the literature for
long memory processes of counts. Hence the question of parametric
estimation has not arisen so far in this context.  However, one may
still be interested in semiparametric estimation of long memory in
counts. We present the following result on the behavior of the
Discrete Fourier Transforms (DFTs) of processes of counts induced by
infinite-variance durations that will be of relevance to us in
understanding the behavior of the GPH estimator. Let $n$ denote the
number of observations on the counts, $\omega_j = 2 \pi j /n$, and
define
\[
J_{n,j}^{\Delta N} = \frac{1}{\sqrt{2\pi n}}\sum_{t'=1}^{n}\Delta
N_{t'} \E^{\I t'\omega_j} \; .
\]
Assume that the distribution of the durations satisfies
\begin{gather} \label{eq:dist_duration}
  P(\tau_k \ge x) \backsim \ell(x)x^{-\alpha} \qquad\mbox{$x
    \rightarrow \infty$}
\end{gather}
where $\ell(x)$ is a slowly varying function with
$\lim_{x\rightarrow\infty}\frac{\ell(kx)}{\ell(x)}=1$ $\forall k>0$ and
$\ell(x)$ is ultimately monotone at $\infty$.

\begin{theorem} \label{theo:stableDFT}
  Let $\{\tau_k\}$ be i.i.d. random variables which satisfy
  (\ref{eq:dist_duration}) with $\alpha \in(1,2)$ and mean
  $\mu_{\tau}$. Then for each fixed $j$, $ \ell(n)^{-1}
  n^{1/2-1/\alpha} J_{n,j}^{\Delta N}$ converges in distribution to a
  complex $\alpha$-stable distribution. Moreover, for each fixed $j$,
  $\omega_j ^d J_{n,j}^{\Delta N} \stackrel{p}{\rightarrow 0}$, where
  $d=1-\alpha/2$.
\end{theorem}

The theorem implies that when $j$ is fixed, the normalized periodogram
of the counts, $\omega_j^{2d} I_{\Delta N}(\omega_j)$ converges in
probability to zero. The degeneracy of the limiting distribution of
the normalized DFTs of the counts suggests that the inclusion of the
very low frequencies may induce negative finite-sample bias in
semiparametric estimators. In addition, the fact that the suitably
normalized DFT has an asymptotic stable distribution could further
degrade the finite-sample behavior of semiparametric estimators, more
so perhaps for the Whittle-likelihood-based estimators than for the
GPH estimator since the latter uses the logarithmic transformation.

By contrast, for linear long-memory processes, the normalized
periodogram has a nondegenerate positive limiting distribution.
See, for example, \cite{terrin:hurvich:1994}.

\subsection{Low-Frequency DFTs of Counts from LMSD Durations}
\label{LowFreqLMSD}

We now study the behavior of the low-frequency DFTs of counts
generated from finite-variance LMSD durations.

\begin{theorem} \label{LMSDDFT}
  Let the durations $\{\tau_k\}$ follow an LMSD model with memory
  parameter $d$. Then for each fixed $j$, $\omega_j^d J_{n,j}^{\Delta
    N}$, converges in distribution to a zero-mean Gaussian random
  variable.

\end{theorem}

This result is identical to what would be obtained if the counts were
a linear long-memory process, and stands in stark contrast to Theorem
\ref{theo:stableDFT}. The discrepancy between these two theorems
suggests that the low frequencies will contribute far more bias to
semiparametric estimates of $d$ based on counts if the counts are
generated by infinite-variance durations than if they were generated
from LMSD durations.

\subsection{Low and High Frequency DFTs of Shot-Noise Processes}
Let $X$ be either the renewal-reward process defined
in~(\ref{eq:defrenewalreward}) or the error duration
process~(\ref{eq:defparkec}). \cite{hsieh:hurvich:soulier:2003},
Theorem 4.1, have proved that Theorem~\ref{theo:stableDFT} still
holds, i.e.  $n^{1/2-1/\alpha} J_{n,j}^{X}$ converges in distribution
to an $\alpha$-stable law, where $\alpha$ is the tail index of the
duration.  This result can probably be extended to all the shot-noise
process for which convergence in distribution of the partial sum
process can be proved.

The DFTs of these processes have an interesting feature, related to
the slow growth/fast growth phenomenon. The high frequency DFTs, i. e.
the DFT $J_{n,j}^{X}$ computed at a frequency $\omega_j$ whose index
$j$ increases as $n^\rho$ for some $\rho>1-1/\alpha$, renormalized by
the square root of the spectral density computed at $\omega_j$, have a
Gaussian weak limit.  This is proved in Theorem 4.2 of
\cite{hsieh:hurvich:soulier:2003}.

\subsection{Estimation of the memory parameter of the LMSV and LMSD models}
We now discuss parametric and semiparametric estimation of the memory
parameter for the LMSV/LMSD models. Note that in both the LMSV and
LMSD models, $\log x_t^2$ can be expressed as the sum of a long memory
signal and $iid$ noise.  Specifically, we have
\begin{equation}\label{eq:SigNoise}
\log X_t^2=\mu +h_t+u_t,
\end{equation}
where $\mu =E\left( \log v_t^2\right) $ and $u_t=\log
v_t^2-E\left( \log v_t^2\right) $ is a zero-mean $iid$ series
independent of $\left\{ h_t\right\} .$ Since all the extant
methodology for estimation for the LMSV model exploits only the
above signal plus noise representation, the methodology continues
to hold for the LMSD model.

Assuming that $\{h_t\}$ is Gaussian, \cite{deo:hurvich:2001} derived
asymptotic theory for the log-periodogram regression estimator (GPH;
\cite{geweke:porter-hudak:1983}) of $d$ based on $\{\log X_t^2\}$.
This provides some justification for the use of GPH for estimating
long memory in volatility. Nevertheless, it can also be seen from
Theorem 1 of \cite{deo:hurvich:2001} that the presence of the noise
term $\{u_t\}$ induces a negative bias in the GPH estimator, which in
turn limits the number $m$ of Fourier frequencies which can be used in
the estimator while still guaranteeing $\sqrt m$-consistency and
asymptotic normality. This upper bound, $m=o[n^{4d/(4d+1)}]$, where
$n$ is the sample size, becomes increasingly stringent as $d$
approaches zero. The results in \cite{deo:hurvich:2001} assume that
$d>0$ and hence rule out valid tests for the presence of long memory
in $\{h_t\}$. Such a test based on the GPH estimator was provided and
justified theoretically by \cite{hurvich:soulier:2002}.

\cite{sun:phillips:2003} proposed a nonlinear log-periodogram
regression estimator $\hat d_{\mathrm{NLP}}$ of $d$, using Fourier
frequencies $1,\ldots,m$. They partially account for the noise term
$\{u_t\}$ through a first-order Taylor expansion about zero of the
spectral density of the observations, $\{\log X_t^2\}$.  They
establish the asymptotic normality of $m^{1/2}(\hat
d_{\mathrm{NLP}}-d)$ under assumptions including $n^{-4d} m^{4d+1/2}
\rightarrow \mathrm{Const}$. Thus, $\hat d_{\mathrm{NLP}}$, with a
variance of order $n^{-4d/(4d+1/2)}$, converges faster than the GPH
estimator, but still arbitrarily slowly if $d$ is sufficiently close
to zero. \cite{sun:phillips:2003} also assumed that the noise and
signal are Gaussian. This rules out most LMSV/LMSD models, since
$\{\log v_t^2\}$ is typically non-Gaussian.

For the LMSV/LMSD model, results analogous to those of
\cite{deo:hurvich:2001} were obtained by \cite{arteche:2004} for the
GSE estimator, based once again on $\{\log X_t^2\}$. The use of GSE
instead of GPH allows the assumption that $\{h_t\}$ is Gaussian to be
weakened to linearity in a Martingale difference sequence.
\cite{arteche:2004} requires the same restriction on $m$ as in
\cite{deo:hurvich:2001}. A test for the presence of long memory in
$\{h_t\}$ based on the GSE estimator was provided by
\cite{hurvich:moulines:soulier:2005}.

\cite{hurvich:ray:2003} proposed a local Whittle estimator of $d$,
based on log squared returns in the LMSV model.  The local Whittle
estimator, which may be viewed as a generalized version of the GSE
estimator, includes an additional term in the Whittle criterion
function to account for the contribution of the noise term
$\{u_t\}$ to the low frequency behavior of the spectral density of
$\{\log X_t^2\}$. The estimator is obtained from numerical
optimization of the criterion function. It was found in the
simulation study of \cite{hurvich:ray:2003}  that the local Whittle
estimator can strongly outperform GPH, especially in terms of bias
when $m$ is large.

Asymptotic properties of the local Whittle estimator were obtained by
\cite{hurvich:moulines:soulier:2005}, who allowed $\{h_t\}$ to be a
long-memory process, linear in a Martingale difference sequence, with
potential nonzero correlation with $\{u_t\}$. Under suitable
regularity conditions on the spectral density of $\{h_t\}$,
\cite{hurvich:moulines:soulier:2005} established the
$\sqrt{m}$-consistency and asymptotic normality of the local Whittle
estimator, under certain conditions on $m$. If we assume that the
short memory component of the spectral density of $\{h_t\}$ is
sufficiently smooth, then their condition on $m$ reduces to
\begin{gather} \label{eq:mclt}
\lim_{n \to \infty} \left(m^{-4d-1+\delta}n^{4d} + n^{-4}
m^{5}\log^2(m) \right) = 0
\end{gather}
for some arbitrarily small $\delta>0$.

The first term in \eqref{eq:mclt} imposes a lower bound on the
allowable value of $m$, requiring that $m$ tend to $\infty$ faster
than $n^{4d/(4d+1)}$. It is interesting that \cite{deo:hurvich:2001},
under similar smoothness assumptions, found that for $m^{1/2} ( \hat
d_{GPH} - d )$ to be asymptotically normal with mean zero, where $\hat
d_{GPH}$ is the GPH estimator, the bandwidth $m$ must tend to $\infty$
at a rate $slower$ than $n^{4d/(4d+1)}$. Thus for any given $d$, the
optimal rate of convergence for the local Whittle estimator is faster
than that for the GPH estimator.

Fully parametric estimation in LMSV/LMSD models once again is
based on $\{\log X_t^2\}$ and exploits the signal plus noise
representation (\ref{eq:SigNoise}). When $\{h_t\}$ and $\{u_t\}$
are independent, the spectral density of $\{\log X_t^2\}$ is
simply the sum of the spectral densities of $\{h_t\}$ and
$\{u_t\}$, viz.
\begin{equation}
f_{\log X^2} (\lambda) = f_{h} (\lambda) + \sigma_u^2 / (2\pi) ,
\end{equation}
where $f_{\log X^2}$ is the spectral density of $\{\log X_t^2\}$,
$f_{h}$ is the spectral density of $\{h_t\}$ and
$\sigma_u^2=\mathrm{var}(u_t)$, all determined by the assumed parametric
model. This representation suggests the possibility of estimating the
model parameters in the frequency domain using the Whittle likelihood.
Indeed, \cite{hosoya:1997} claims that the resulting estimator is
$\sqrt n$-consistent and asymptotically normal. We believe that though
the result provided in \cite{hosoya:1997} is correct, the proof is
flawed. \cite{deo:1995o} has shown that the quasi-maximum likelihood
estimator obtained by maximizing the Gaussian likelihood of $\{\log
X_t^2\}$ in the time domain is $\sqrt n$-consistent and asymptotically
normal.

One drawback of the latent-variable LMSV/LMSD models is that it is
difficult to derive the optimal predictor of $|X_t|^s$. In the LMSV
model, $\{|X_t|^s\}$ for $s>0$ serves as a proxy for volatility, while
in the LMSD model, $\{X_t\}$ represents durations. A computationally
efficient algorithm for optimal linear prediction of such series was
proposed in \cite{deo:hurvich:lu:2005}, exploiting the Preconditioned
Conjugate Gradient (PCG) algorithm. In \cite{chen:hurvich:lu:2005}, it
is shown that the computational cost of this algorithm is $O(n
\log^{5/2} n )$, in contrast to the much more expensive Levinson
algorithm, which has cost of $O(n^2)$.

\subsection{Simulations on the GPH Estimator for Counts}

\begin{figure}
\caption{Averaged $\log-\log$ periodogram plots for the counts generated
  from $iid$ Stable and LMSD durations.}
\begin{center}
  \includegraphics[width=8cm ] {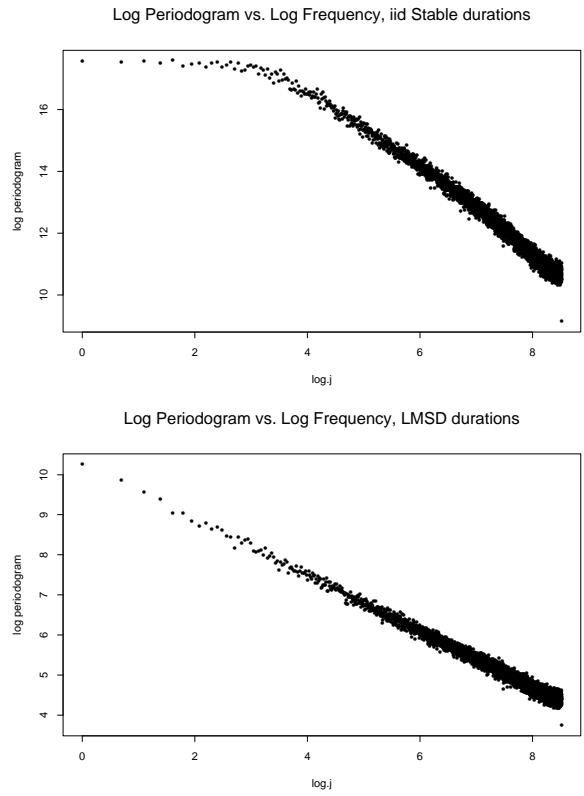}
\end{center}
\label{fig:log_log_pgram}
\end{figure}

We simulated i.i.d. durations from a positive stable distribution
with tail index $\alpha=1.5$, with an implied $d$ for the counts
of $.25$. We also simulated durations from an LMSD $(1,d,0)$ model
with Weibull innovations, $AR(1)$ parameter of $-.42$, and
$d=.3545$, as was estimated from actual tick-by-tick durations in
\cite{deo:hsieh:hurvich:2005}. The stable durations were
multiplied by a constant $c=1.21$ so that the mean duration
matches that found in actual data. For the LMSD durations, we used
$c=1$. One unit in the rescaled durations is taken to represent
one second. Tables $\ref{tab:SCD_0215}$ and
$\ref{tab:LMSD_sim_count}$, for the stable and LMSD cases
respectively, present the GPH estimates based on the resulting
counts for different values of $\Delta t$, using $n=10,000$,
$m=n^{0.5}$ and $m=n^{0.8}$. For the stable case, the bias was far
more strongly negative for the smaller value of $m$, whereas for
the LMSD case, the bias did not change dramatically with $m$. This
is consistent with the discussion in Section \ref{LowFreqLMSD},
and also with the averaged $\log-\log$ periodogram plots presented
in Figure \ref{fig:log_log_pgram}, where the averaging is taken
over a large number of replications, and all positive Fourier
frequencies are considered, $j=1,\ldots,n/2$. The plot for the
stable durations (upper panel) shows a flat slope at the low
frequencies. For this process, using more frequencies in the
regression seems to mitigate the negative bias induced by the
flatness in the lower frequencies as indicated by the less biased
estimates of $d$ when $m=n^{0.8}$.

For the LMSD process, if the conjecture is correct then the counts
should have the same memory parameter as the durations, $d=.3545$.
Assuming that this is the case, we did not find severe negative
bias in the GPH estimators on the counts, though the estimate of
$d$ seems to increase with $\Delta t$ in the case when
$m=n^{0.5}$. The averaged $\log-\log$ periodogram plot presented
in the lower panel of Figure \ref{fig:log_log_pgram} shows a
near-perfect straight line across all frequencies, which is quite
different from the pattern we observed in the case of counts based
on stable durations. The straight-line relationship here is
consistent with the bias results in our LMSD simulations, and with
the discussion in Section \ref{LowFreqLMSD}.

Statistical properties of $\hat{d}_{GPH}$ and the choice of $m$
for Gaussian long-memory time series have been discussed in recent
literature. \cite{robinson:1995l} showed for Gaussian processes that
the GPH estimator is $m^{1/2}$-consistent and asymptotically
normal if an increasing number of low frequencies $L$ is trimmed
from the regression of the log periodogram on log frequency.
\cite{hurvich:deo:brodsky:1998} showed that trimming can be
avoided for Gaussian processes. In our simulations, we did not use
any trimming. There is as yet no theoretical justification for the
GPH estimator in the current context since the counts are clearly
non-Gaussian, and presumably constitute a nonlinear process. It is
not clear whether trimming would be required for such a theory,
but our simulations and theoretical results suggest that in some
situations trimming may be helpful, while in others it may not be
needed.

%\clearpage

\begin{table}
\caption{GPH estimators for counts with different $\triangle t$.
Counts generated from $iid$ stable durations with skewness
parameter $\beta=0.8$ and tail index $\alpha=1.5$. The
corresponding memory parameter for counts is $d=.25$. We generated
500 replications each with sample size $n=10,000$. The number of
frequencies in the log periodogram regression was $m=n^{0.8}=1585$
and $m=\sqrt{n}=100$. $t$-values marked with $\ast$ reject the
null hypothesis, $d = 0.25$ in favor of $d<0.25$.}
\begin{center}
\label{tab:SCD_0215}
\begin{tabular}{|c|c|c|c|c|c|c|}\hline
$\triangle t$ & \multicolumn{2}{c|}{$m=n^{0.5}$}&
\multicolumn{2}{c|}{$m=n^{0.8}$}
\\[0.5ex]\hline
$c=1.21$ & $Mean(\widehat d_{GPH})$ & t-Value & $Mean(\widehat
d_{GPH})$ & t-Value
\\[0.5ex]
\hline\hline
 5 min & 0.1059 & $-17.65^{\ast}$ & 0.2328 & $-5.77^{\ast}$ \\
10 min & 0.0744 & $-23.08^{\ast}$ & 0.2212 & $-8.31^{\ast}$ \\
20 min & 0.0715 & $-23.23^{\ast}$ & 0.2186 & $-7.75^{\ast}$ \\
\hline
\end{tabular}
\end{center}
\end{table}

\begin{table}
\caption{Mean of the GPH estimators for counts with different
$\Delta t$. Counts generated from LMSD durations with Weibull (1,
$\gamma$) shocks. The number of frequencies in the log periodogram
regression was $m=\sqrt{n}$ and $m=n^{0.8}$. We used $d=.3545$ and
$\gamma=1.3376$ for our simulations. We simulated $200$
replications of the counts, each with sample size $n=10,000$.
$t$-values marked with $\ast$ reject the null hypothesis,
$d=0.3545$ in favor of $d < 0.3545$.}{\label{tab:LMSD_sim_count}}
\begin{center}
\begin{tabular}{|c|c|c|c|c|c|c|}\hline
$\triangle t$ & \multicolumn{2}{c|}{$m=n^{0.5}$}&
\multicolumn{2}{c|}{$m=n^{0.8}$}
\\[0.5ex]\hline
$c=1$ & $Mean(\widehat d_{GPH})$ & t-Value & $Mean(\widehat
d_{GPH})$ & t-Value
\\[0.5ex]
\hline\hline
 5 min & 0.3458 & $-1.76^{\ast}$ & 0.3471 & $-6.49^{\ast}$ \\
30 min & 0.3873 & $3.45^{\ast}$ & 0.3469 & $-3.59^{\ast}$ \\
60 min & 0.3923 & $4.05^{\ast}$ & 0.3478 & $-3.20^{\ast}$ \\
\hline
\end{tabular}
\end{center}
\end{table}

%\clearpage
\clearpage

\subsection{Estimation of the memory parameter 
  of the Infinite Source Poisson process} Due to the underlying
Poisson point process, the Infinite Poisson Source process is a very
mathematically tractable model. Computations are very easy and in
particular, convenient formulas for cumulants of integrals along paths
of the process are available. This allows to derive the theoretical
properties of estimators of the Hurst index or memory parameter.
\cite{fay:roueff:soulier:2005} have defined an estimator of the Hurst
index of the Infinite Poisson source process (with random transmission
rate) related to the GSE and proved its consistency and rate of
convergence.  Instead of using the DFTs of the process, so-called
wavelets coefficients are defined as follows. Let $\psi$ be a
measurable compactly supported function on $\mathbb{R}$ such that $\int
\psi(s) \, \D s = 0$.  For $j\in\mathbb{N}$ and $k=0,\dots,2^j-1$,
define
\begin{gather*}
  w_{j,k} = \int \phi(s) X_s \, \D s \; .
\end{gather*}
If~(\ref{eq:nonindep}) holds, then $\mathbb{E}[w_{j,k}] = 0$ and
$\mathrm{var}(w_{j,k}) = L(2^j) 2^{(2-\alpha)j} = L(2^j) 2^{2dj}$,
where $\alpha$ is the tail index of the durations, $d=1-\alpha/2$ is
the memory parameter and $L$ is a slowly varying function at infinity.
This scaling property makes it natural to define a contrast function
\begin{gather*}
  \hat W(d') = \log \left( \sum\nolimits_{(j,k)\in\Delta} 2^{-2d'j} w^{2}_{j,k}
  \right) + \delta d' \log(2) \; ,
\end{gather*}
where $\Delta$ is the admissible set of coefficients, which depends on
the interval of observation and the support of the function $\psi$.
The estimator of $d$ is then $\hat d = \arg\min_{d'\in(0,1/2)} W(d')$.
\cite{fay:roueff:soulier:2005} have proved under some additional
technical assumptions that this estimator is consistent. The rate of
convergence can be obtained, but the asymptotic distribution is not
known, though it is conjectured to be Gaussian, if the set $\Delta$ is
properly chosen. 

Note in passing that here again, the slow growth/fast growth
phenomenon arises. It can be shown, if the shocks and durations are
independent, that for fixed $k$, $2^{(1-\alpha)j/2} w_{j,k}$ converges
to an $\alpha$-stable distribution, but if $k$ tends to infinity at a
suitable rate, $2^{-dj} w_{jk}$ converges to a complex Gaussian
distribution. This slow growth/fast growth phenomenon is certainly a
very deep property of these processes that should be understood more
deeply.

%%%%%%%%%%%%%%%%%%%%%%%%%%%%%%%%%%%%%%%%%%%%%%%%%%%%%%%%%%%%%%%%%%%%%%

\section*{Appendix}
\begin{proof}[ of Theorem \ref{theo:stableDFT}]
  For simplicity, we set the clock-time spacing $\Delta t=1$.  Define
\begin{gather*}
  S_{\tau,n}(\theta)=\sum_{k=1}^{\lfloor n\theta \rfloor}\tau_k
  \qquad 0 \leq \theta \leq 1 \; , \\
  S_{\Delta N,n}(\theta)=\sum_{t'=1}^{\lfloor n\theta \rfloor} \Delta
  N_{t'} \qquad 0 \leq \theta \leq 1 \; .
\end{gather*}
Since $\alpha <2$ and $\{\tau_k\}$ is an i.i.d. sequence, by the
fonctional central limit theorem (FCLT) for random variables in the
domain of attraction of a stable law (see \cite[Theorem
2.4.10]{embrechts:kluppelberg:mikosch:1997}), $l(n) n^{-1/\alpha}
\{ S_{\tau,n}(\theta)-\lfloor n \theta \rfloor \mu_{\tau}\}$ converges
weakly in $\mathcal{D}(0,1)$ to an $\alpha$-stable motion, for some
slowly varying function $l$. Now define
\[
U_n(\theta) = (2\pi)^{-1/2} l(n) n^{-1/\alpha} \{S_{\Delta N,n}(\theta)-\lfloor
n\theta \rfloor/\mu_{\tau}\} \; .
\]
By the equivalence of FCLTs for the counting process and its
associated partial sums of duration process (see
\cite{iglehart:whitt:1971}), $U_n$ also converges weakly in
$\mathcal{D}([0,1])$ to an $\alpha$-stable motion, say $S$.  Summation
by parts yields, for any nonzero Fourier frequency $\omega_j$ (with
fixed $j>0$)
\begin{align*}
  l(n) & n^{1/2-1/\alpha} J^{\Delta N}_{n,j} = (2\pi)^{-1/2} l(n)
  n^{-1/\alpha}\sum_{t'=1}^{n} \{\Delta N_{t'} -
  1/\mu_\tau\} \, \E^{\I t' \omega_j}  \\
  & = \sum_{t'=1}^{n} \{U_n(t'/n) - U_n((t'-1)/n)\} \, \E^{\I
    t'\omega_j}
%%  - l(n) n^{-1/\alpha}\Delta N_n  \E^{\I n \omega_j} \\  & 
= \int_0^1 \E^{2 \I j \pi x} \, \D U_n(x) 
  %%  \ - \ l(n) n^{-1/\alpha}\Delta N_n \E^{\I n \omega_j}
\; .
\end{align*}
%% Note that $ l(n) n^{-1/\alpha} \Delta N_n = o_p(1)$. 
Hence by the continuous mapping theorem
\[
\sqrt{2\pi } \; l(n) n^{1/2-1/\alpha} J^{\Delta N}_{n,j}
\stackrel{d}{\longrightarrow} \int^1_0 \E^{2 \I \pi j x} \, \D S(x)
\]
which is a stochastic integral with respect to a stable motion, hence
has a stable law.

To prove the second statement of the theorem, note that for fixed $j$
and as $n \rightarrow \infty$, $f(\omega_j) \sim l_1(n)
\omega^{-2d}_j$ for some slowly varying function $l_1$, so
\begin{multline}\label{eq:norm_pgram}
  f^{-1/2}(\omega_j)J^{\Delta N}_{n,j} =
  \frac{l(n)n^{1/\alpha-1/2}}{f^{1/2}(\omega_j)}\frac{J^{\Delta
      N}_{n,j}}{l(n)n^{1/\alpha-1/2}}  \\
  \sim C_1l(n)n^{1/\alpha+\alpha/2-3/2}\frac{J^{\Delta
      N}_{n,j}}{\mu_{\tau}^{-1-1/\alpha}l(n)n^{1/\alpha-1/2}} \; .
\end{multline}
Since $1/\alpha+\alpha/2-3/2<0$, we have
$l(n)n^{1/\alpha+\alpha/2-3/2}\rightarrow 0$. Hence by Slutsky's
Theorem, (\ref{eq:norm_pgram}) converges to zero. \qed
\end{proof}

\begin{proof} [of Theorem \ref{LMSDDFT}]
  Let $S_n(t)=n^{-H}\sum_{k=1}^{[nt]} (\tau_k-\mathbb{E}[\tau_k])$, $t
  \in (0,1)$. It is shown in Surgailis and Viano (2002) that $S_n(t)
  \stackrel{d}{\Rightarrow} B_{H}(t)$ in $\mathcal{D}([0,1])$ where
  $B_{H}(t)$ is fractional Brownian motion with Hurst parameter
  $H=d+1/2$. Thus, by Iglehart and Whitt (1971), it follows that
  $t^{-H}N \rightarrow A B_{H}$ in $\mathcal{D}([0,1])$, where $A$ is
  a nonzero constant. The result follows as above by the continuous
  mapping theorem and summation by parts.  \qed
\end{proof}

\end{document}